\documentclass[a4paper,11pt,oneside,reqno]{amsart}
\usepackage[all,cmtip]{xy}
\usepackage{amsmath}
\usepackage{amsthm}
\usepackage{amssymb}
\usepackage{graphicx}
\usepackage{mathrsfs}
\usepackage{bbm}
\theoremstyle{plain}
\newtheorem{teore}{Theorem}[section]
\newtheorem{defin}{Definition}[section]
\newtheorem{lem}[teore]{Lemma}
\newtheorem{coro}[teore]{Corollary}
\newtheorem{propo}[teore]{Proposition}

\newtheorem{claim}{Claim}[teore]
\newtheorem*{claim*}{Claim}
\newtheorem*{thm*}{Theorem}
\newtheorem*{defi*}{Definition}
\theoremstyle{remark}
\newtheorem{ejemplo}[teore]{{\sc Example}}
\newtheorem{ejemplos}[teore]{{\sc Examples}}
\newtheorem{notas}[teore]{{\sc Remark}}
\newcommand{\nrm}[1]{\|#1\|}

\newcommand{\prop}{\begin{propo}}
\newcommand{\fprop}{\end{propo}}
\newcommand{\cor}{\begin{coro}}
\newcommand{\fcor}{\end{coro}}

\newcommand{\defi}{\begin{defin}}
\newcommand{\fdefi}{\end{defin}}
\newcommand{\eje}{\begin{ejemplo}}
\newcommand{\feje}{\end{ejemplo}}
\newcommand{\ejes}{\begin{ejemplos}}
\newcommand{\fejes}{\end{ejemplos}}
\newcommand{\lema}{\begin{lem}}
\newcommand{\flema}{\end{lem}}
\newcommand{\teor}{\begin{teore}}
\newcommand{\fteor}{\end{teore}}
\newcommand{\nota}{\begin{notas}}
\newcommand{\fnota}{ \end{notas}}
\newcommand{\clam}{\begin{claim}}
\newcommand{\fclam}{\end{claim}}
\newcommand{\clams}{\begin{claim*}}
\newcommand{\fclams}{\end{claim*}}

\newcommand{\lclam}{\begin{lclaim}}
\newcommand{\flclam}{\end{lclaim}}
\newcommand{\prucl}{\prue[Proof of Claim:]}
\newcommand{\fprucl}{\fprue}
\newcommand{\ben}{\begin{enumerate}}
\newcommand{\een}{\end{enumerate}}
\newcommand{\bit}{\begin{itemize}}

\newcommand{\eit}{\end{itemize}}

\newcommand{\mc}[1]{\mathcal{#1}}

\newcommand{\mr}[1]{\mathrm{#1}}

\newcommand{\mk}[1]{\mathfrak{#1}}

\newcommand{\we}{\wedge}

\renewcommand{\bar}[1]{\overline{#1}}

\newcommand{\casos}{\begin{itemize}}
\newcommand{\fcasos}{\end{itemize}\setcounter{cs}{1}}

\newcommand{\conj}[2]{ \{ {#1}\,:\,{#2} \} }

\newcommand{\om}{\omega}

\newcommand{\buit}{\emptyset}

\newcommand{\ga}{\gamma}

\newcommand{\al}{\alpha}

\newcommand{\de}{\delta}
\newcommand{\De}{\Delta}

\newcommand{\vep}{\varepsilon}
\renewcommand{\S}{\mathcal{S}}

\newcommand{\N}{{\mathbb N}}

\newcommand{\rest}{\upharpoonright}

\newcommand{\con}{\subseteq}

\newcommand{\prue}{\begin{proof}}
\newcommand{\fprue}{\end{proof}}
\begin{document}

\setlength{\baselineskip}{6mm}

\title{Shellable weakly compact subsets of $C[0,1]$}

\author[J. Lopez-Abad]{J. Lopez-Abad}
\address{J. Lopez-Abad \\ Instituto de Ciencias Matematicas (ICMAT).  CSIC-UAM-UC3M-UCM. C/ Nicol\'{a}s Cabrera 13-15, Campus Cantoblanco, UAM
28049 Madrid, Spain} \email{abad@icmat.es}

\author[P. Tradacete]{P. Tradacete}
\address{P. Tradacete\\Mathematics Department\\ Universidad Carlos III de Madrid \\  28911 Legan\'es (Madrid). Spain.}
\email{  ptradace@math.uc3m.es }

\thanks{Support of the Ministerio de Econom\'{\i}a y Competitividad under grant MTM2012-31286 (Spain) is gratefully acknowledged. First author has been partially supported by ICMAT Severo Ochoa project SEV-2015-0554 (MINECO).
 Second author also supported by grants MTM2013-40985 and Grupo UCM 910346.  }



\begin{abstract}
We show that for every weakly compact subset $K$ of $C[0,1]$ with finite Cantor-Bendixson rank, there is a reflexive Banach lattice $E$ and an operator $T:E\rightarrow C[0,1]$ such that $K\subseteq T(B_E)$. On the other hand, we exhibit an example of a weakly compact set of $C[0,1]$ homeomorphic to $\omega^\omega+1$ for which such $T$ and $E$ cannot exist. This answers a question of M. Talagrand in the 80's.
\end{abstract}

\maketitle

\section{Introduction}

In the celebrated paper \cite{DFJP}, W. Davis, T. Figiel, W. Johnson and A. Pelczynski showed that every weakly compact operator between Banach spaces factors through a reflexive Banach space. This fact relies on a general construction, arising from interpolation theory, which shows that every weakly compact set in a Banach space is contained in the image of the unit ball by an operator defined on a reflexive Banach space. Analogous statements in the category of Banach lattices and Banach algebras have been considered respectively by C. Aliprantis and O. Burkinshaw in \cite{AB:84} and A. Blanco, S. Kaijser and T. J. Ransford in \cite{BKR}. In the former paper, it was shown that under some extra assumptions, a weakly compact operator between Banach lattices can be factored through a reflexive Banach lattice. In \cite{T}, M. Talagrand showed that these extra assumptions cannot be completely removed.

In fact, let us say that a set $A$ in a Banach space $X$ is shellable by reflexive Banach lattices if there
is a reflexive Banach lattice $E$ and an operator $T:E\rightarrow X$ such that $A\subseteq T(B_E)$, where
$B_E$ denotes the unit ball of $E$. In \cite{T}, the author constructs a weakly compact set $K_{\mathcal T}$
of continuous functions on the interval $[0,1]$ which is not shellable by reflexive Banach lattices.
Associated to this compact set $K_{\mathcal T}$, one can consider a weakly compact operator
$T:\ell_1\rightarrow C[0,1]$ which cannot be factored through a reflexive Banach lattice. The compact set
$K_{\mathcal T}$ is small, as it is homeomorphic to the ordinal $\om^{\om^2}+1$. It was naturally  asked in
\cite{T} what  is the smallest ordinal $\alpha$ for which there is a weakly compact set $K\subset C[0,1]$
homeomorphic to $\alpha$ and not shellable by  reflexive Banach lattices. The main result of this note is
that this ordinal is precisely $\omega^\omega+1$.

Our proof has  two parts. For the lower bound, we will use crucially that when $K$ is a countable compact
space with finite Cantor-Bendixson rank then the space $C(K)$ of continuous functions on $K$ has an
unconditional basis, in fact it is isomorphic to $c_0$. This implies that when $K$ is a countable
weakly-compact subset of $C[0,1]$,   the natural evaluation  operator $\xi\in C[0,1]^* \mapsto  C(K)$,
$\xi(\mu)(p):= \int p d\mu$, $p\in K$, factors through a  reflexive Banach \emph{lattice}. From this, a simple
duality argument then shows that $K$ is shellable  by a reflexive Banach lattice.

For the upper bound we will construct a weakly compact set $K_\omega\subset C(\om^{\om^2}+1)$  homeomorphic
to $\omega^\omega$  which is not shellable by reflexive Banach lattices. It is well-known that on $L_1(\mu)$ spaces, weak-compactness, uniform integrability and the Banach-Saks property are equivalent. In \cite{T}, uniform integrability in combination with martingale techniques was heavily used. In this paper, instead of that we exploit the weak Banach-Saks property of $L_1$.

\section{Preliminaries}
We shall make use of standard Banach space facts and terminology as may be
found in \cite{AK, LT1, LT2}. Throughout, $B_X$ denotes the unit ball of a Banach space $X$, that is $B_X=\{x\in X:\|x\|\leq1\}$. Also, an operator will always refer to a bounded and linear map. Recall that $A\subset X$ is a Banach-Saks set if every $(x_n)\subset A$ has a Cesaro convergent subsequence $(x_{n_k})$, that is the arithmetic means $\frac1n\sum_{k=1}^nx_{n_k}$ are convergent. Recall that a reflexive Banach lattice $E$ can
always be represented as a space of measurable functions on some probability space $(\Omega,\Sigma,\mu)$ such
that $L^\infty(\mu)\subset E\subset L^1(\mu)$, with $\|f\|_1\leq \|f\|_E$ \cite[Theorem 1.b.12]{LT2}. In particular, $E^*$ is
also a space of measurable functions on $(\Omega,\Sigma,\mu)$ and both $B_E$ and $B_{E^*}$ are
equi-integrable sets in $L^1(\mu)$ (or equivalently, Banach-Saks sets). A subset of a Banach lattice
$A\subseteq E$ is solid when $|x|\le|y|$ with $y\in A$ implies that $x\in A$. Given a set $A\subseteq E$, we
define its solid hull $\mr{sol}(A)=\conj{x\in E}{\,|x|\leq |y|  \text{ for some $y\in A$}}$.  The reader is
refered to \cite{AB, LT2} for further explanations concerning Banach lattices.

Recall the classical result by S. Mazurkiewicz and W. Sierpinski \cite{MS} that states that every countable compact
space $K$ is homeomorphic to a unique ordinal number $\om^{\al_K} \cdot n_k +1$, $n_K\in \N$, with  its order
topology.  We will use here the fact that when $K$ has finite \emph{Cantor-Bendixson} index, that is when
$\al_K$ is finite, then the corresponding space of continuous functions $C(K)$  has an unconditional basis,
in fact, $C(K)$ is isomorphic to $c_0$. This is actually a characterization of those compact spaces (see \cite[Theorem 4.5.2]{AK} for
more details).

\defi
Let $Y$ be a Banach space, and $\mathcal X$ a collection of Banach spaces. A set $K\subset Y$ will be called $\mc X$-shellable (or shellable by $\mc X$) if there is a space $X$ in $\mc X$ and an operator $T:X\rightarrow Y$ such that $K\subset T(B_X)$. In this case, we will write $K\in \mr{Sh}_X(Y)$.
\fdefi

We will denote by $\mr{Sh}_{\mc X}(Y)$ the family of subsets $K\subset Y$ which are $\mc X$-shellable.
In particular, let $\mr{Sh}_{\mc{R}}(Y)$ and $\mr{Sh}_{\mc{RL}}(Y)$ be the collections of subsets of $Y$ which are shellable by a reflexive Banach space and by a reflexive Banach lattice, respectively. The following is easy to prove.

\prop
\begin{enumerate}
\item[(a)] $\mr{Sh}_{\mc X}(Y)$  is closed under  subsets and   convex hulls.
 \item[(b)] If $\mc X$ is closed under finite direct sums, then $\mr{Sh}_\mc X(Y)$ is closed under finite unions.
 \item[(c)] If $Y_0$ is a complemented subspace of $Y_1$, then $\mr{Sh}_X(Y_0)=\conj{W\con Y_0}{W\in \mr{Sh}_X(Y_1)}$.  \qed
\end{enumerate}
\fprop

We recall the main construction in \cite{DFJP},    later adapted in \cite{AB:84} to the lattice setting (see also \cite[Section 5.2]{AB}).

\teor\label{reflexive-interpolation}
Let $W$ be a convex, circled, bounded subset of a Banach space $X$. There is Banach space $Y$, and a continuous linear injection $J:Y\rightarrow X$ such that $W\subset J(B_Y)$. In addition,
\begin{enumerate}
\item[(a)] $Y$ is reflexive if and only if $W$ is relatively weakly compact.
\item[(b)] If $X$ is a Banach lattice, and $W$ is solid, then $Y$ can be taken to be a Banach lattice.
 \end{enumerate}

\fteor
It is not true in general that the solid hull of a weakly-compact set in a Banach lattice is weakly-compact. Nevertheless, if $X$ is a space with an unconditional basis, or if $X$ is a Banach lattice not containing $c_0$, then every weakly compact set has weakly compact solid hull (see \cite{CW} for more details).  As a consequence we obtain the following.

\cor
Let $W$ be a relatively weakly-compact subset of a Banach space $X$. Then
\begin{enumerate}
\item[(a)]  $W$ is shellable by a reflexive Banach space.
\item[(b)]  If $X$ is a Banach lattice and $W$ is solid, then $W$ is shellable by a reflexive Banach lattice.
\item[(c)] If $X$  has an  unconditional basis, then $W$ is shellable by a reflexive Banach lattice.
\item[(d)]  If $X$ is a Banach lattice that does not contain $c_0$, then $W$ is shellable by a reflexive Banach lattice. \qed
 \end{enumerate}
\fcor

\section{Main results}

\subsection{Lower bound}
Let us begin this section with the lower bound concerning Talagrand's problem:

\teor\label{kfksdfsklfsd}
 Every countable weakly-compact subset of $C[0,1]$ of finite Cantor-Bendixson rank is in $\mr{Sh}_{\mc RL} (C[0,1])$.
\fteor

\prue
Let us consider a countable weakly-compact subset $K\con C[0,1]$ with finite Cantor-Bendixson rank. Let
$$
\phi:C[0,1]^*\rightarrow C(K)
$$
given by $\phi(\mu) (f)=\int fd\mu$ for each $\mu\in C[0,1]^*$ and $f\in K\subset C[0,1]$. Clearly, $\phi$ defines a bounded linear operator with $\|\phi\|\leq\max\{\|f\|:f\in K\}.$ 

Let $L=\phi(B_{C[0,1]^*})$, which is a weakly compact subset of $C(K)$. Since $K$ has finite Cantor-Bendixson rank, $C(K)$ has an unconditional basis $(e_n)$. Hence, by \cite[Theorem 4.39]{AB}, it follows that its solid hull (with respect to the unconditional basis)
$$
\mr{sol}(L)=\conj{\sum_n a_n e_n\in  C(K)}{ \text{there is $\sum_n b_n e_n\in L$ such that $|b_n|\leq |a_n|$ for all $n$}}
$$
is a weakly compact subset of $C(K)$. Hence, by \cite[Theorem 2.2]{AB:84}, $\phi$ factors through a reflexive Banach lattice $E$. More precisely, there is $\widetilde{\phi}:C[0,1]^*\rightarrow E$ and an interval preserving lattice homomorphism  $J:E\rightarrow C(K)$ (with respect to the lattice structure induced by the unconditional basis $(e_n)$ of $C(K)$) with $\mr{sol}(L)\subset J(B_E)$  and $\|J\|\leq \|\phi\|$, such that the following diagram commutes:
$$
\xymatrix{C[0,1]^*\ar_{\widetilde{\phi}}[dr]\ar[rr]^\phi&&C(K)\\
&E\ar_{J}[ur]& }
$$
Thus, we clearly have the dual diagram as follows:
$$
\xymatrix{C(K)^*\ar_{J^*}[dr]\ar[rr]^{\phi^*}&&C[0,1]^{**}\\
&E^*\ar_{\widetilde{\phi}^*}[ur]& }
$$
Let $j:C[0,1]\rightarrow C[0,1]^{**}$ denote the canonical embedding in the bidual. Note that for every $f\in K$ we have that $\phi^*(\delta_f)=j(f)$. Indeed, for every $\mu\in C[0,1]^*$ we have
$$
\phi^*(\delta_f)(\mu)= \delta_f(\phi(\mu))=\mu(f)=j(f)(\mu).
$$
Now, since $\phi^*$ is $w^*-w$ continuous, we have that
$$
\phi^*(B_{C(K)^*})=\phi^*\big(\overline{\rm co}^{w^*}(\{\pm\delta_f:f\in K\})\big)\subset \overline{\rm co}^{w}(\{\pm\phi^*(\delta_f):f\in K\})\subset j(C[0,1]),
$$
so we actually have $\phi^*:C(K)^*\rightarrow C[0,1]$. Let $F$ be the (reflexive) sublattice of $E^*$ generated by $\{J^*(\delta_f):f\in K\}.$ Let $T:F\rightarrow C[0,1]$, be the restriction of $\widetilde{\phi}^*$ on $F$. Since for $f\in K$ we have that $J^*(\delta_f)\in \|\phi\| B_{F}$, and
$$
T(J^*(\delta_f))=\phi^*(\delta_f)=f,
$$
we have that $K\subseteq \|\phi\|T(B_{F})$ and the proof is finished.
\fprue

\subsection{Upper bound}
We continue by proving that the Cantor-Bendixson rank $\om+1$ is sharp. The main result is the following.
\teor\label{lksfkkwelewre}
There is a weakly-compact subset of $C[0,1]$ homeomorphic to $\om^\om+1$ which is not shellable by reflexive Banach
lattices.
\fteor
In fact, we find a weakly compact subset of $C(\om^{\om^2}+1)$ homeomorphic to $\om^\om+1$ which is not
shellable by reflexive  lattices. This implies Theorem \ref{lksfkkwelewre}, since for any countable compact
space $L$ its space of continuous functions $C(L)$ can be isometrically embedded in $C[0,1]$ in a
complemented  way (any such compact space is homeomorphic to a retract of $\{0,1\}^\N$).   It is more
convenient to work with the appropriate compact families of finite subsets of $\N$ rather than directly with
ordinal numbers. Recall that a family $\mc F$ of finite subsets of an index set $I$ is considered as a topological space
by identifying each element $s$ of $\mc F$ with its characteristic function, and then by considering the
induced product topology on $2^I$. In particular, we will consider the \emph{Schreier} family
$$
\mc S=\{s\subset \mathbb N:\sharp s\leq \min s\}.
$$
It is well-known that $\mc S$ is   compact and  homeomorphic to $\omega^\omega+1$. The Schreier barrier
$\mathfrak{S}$ is the family of maximal elements in $\mc S$ or equivalently, the set of isolated points of
$\mc S$. Explicitly,
$$
\mathfrak{S}=\{s\subset \mathbb N:\sharp s=\min s\}.
$$

Since $\mc S$ is a scattered compact space, its set of isolated points $\mk S$ is dense in it. So it makes sense to consider the following continuous extension property. Let us call a sequence $(x_s)_{s\in \mc S}$ a \emph{weakly-convergent tree} if the assignment $s\in \mc S\mapsto x_s\in E$ is continuous with respect to the weak topology in $E$.

\lema\label{weak tree}
Given a reflexive Banach space $E$ and $(x_s)_{s\in \mk S}$ bounded in $E$, there is an infinite set $M\con \N$  such that we can extend $(x_s)_{s\in \mc S\rest M}$ to be a weakly convergent tree.
\flema

\prue
First of all, we can assume that $(E,w)$ is metrizable with $d$. An induction argument on $n\in \N$ gives that for every $\varepsilon>0$, $n\in \N$ and any infinite set $M\con \N$, there exists an infinite set $N\con M$ such that $(x_s)_{s\in [N]^{\leq n}}$ is a weakly convergent tree with
$$
\mr{diam}_d((x_s)_{s\in [N]^{\leq n}})<\varepsilon.
$$
Using this, one can find a sequence $(M_n)_{n\in N}$ with $M_{n+1}\con M_n\con \N$, so that if $m_n=\min M_n$, then $(x_{\{m_n\}\cup s})_{s\in [M_n]^{\leq m_n-1}}$ is a weakly convergent tree with $d$-diameter smaller than $1/n$. Let $L\con \{m_n\}_{n\in\N}$ be such that $x_{m}\rightarrow_{m\in L} x_\emptyset$. Then $(x_s)_{s\in \mc S\rest L}$ is a weakly convergent tree.
\fprue

The previous Lemma generalizes to uniform families with a similar proof (see \cite{AT} for details on uniform families).

We will also use the ``square'' of the Schreier family $\mc S_2$. Recall that given two families $\mc F$ and $\mc G$ of finite subsets of $\N$ their product is $$\mc F\otimes \mc G:=\conj{\bigcup_{i=1}^n s_i}{s_1<\cdots<s_n,\, s_i\in \mc F\textrm{ for }1\leq i\leq n \text{ and $\{\min s_i\}_i \in \mc G$}
 },
$$
where $s<t$ means $\max s <\min t$.
Let $\mc S_2:= \mc S\otimes \mc S$.  Given $n\in \N$, let $\N_n:=[n,\infty[$. Let $\mc F_n:= [\N_n]^{\le n}$, $\mc G_n:=\mc S\otimes [\N_n]^{\le n}$, $\mc F_\om := \mc S$ and $\mc G_\om:= \mc S_2$. Notice that $\bigcup_n\mc F_n=\mc F_\om$ and $\bigcup_n\mc G_n= \mc G_\om$.
 The families $\mc F_\al$, $\mc G_\al$, $\al\le \om$ are compact, hereditary and any restriction of them have rank $\al+1$ and $\om\cdot \al+1$, respectively.   It is easy to show  from the definition of $\mc S_2$ that
each element $t\in \mc S_2=\mc S\otimes \mc S$ has a unique decomposition $s=s[0]\cup s[1]\cdots \cup s[n]$, where\begin{enumerate}
\item[(a)] $s[0]<s[1]<\cdots <s[n]$,
\item[(b)] $\{\min s[i]\}_{i\le n}\in \mc S$,
\item[(c)] $s[0],\cdots,s[n-1]\in\mathfrak S$ and $s[n]\in \mc S$.
\end{enumerate}

\defi
Given $s=\{m_0<\cdots <m_k\}\in \mc S$ and $t=t[0]\cup \cdots \cup t[l]\in \mc S_2$, let us denote
$$
\langle s,t\rangle=\#(\conj{0\le i \le \min\{k,l\}}{ m_i \in t[i]}).
$$
Let $\Theta: \mc S \times \mc S_2 \to \{0,1\}$ be the mapping that to $(s,t)\in \mc S\times \mc S_2$ assigns
 \begin{align}
 \Theta(s,t):= & \langle s,t\rangle +1 \qquad \mod 2.
 \end{align}
 \fdefi

\prop\label{kjndfksdfeww}
The mapping $\Theta$ is coordinate-wise continuous.
\fprop
\prue
Since
$$\Theta(s,t)=    \frac12\Big((-1)^{\langle s,t\rangle} +1\Big),$$
we only need to prove that   $\langle \cdot,\cdot \rangle$ is a coordinate-wise continuous mapping.
Suppose that $s_n\to_n s$ in $\mc S$, and fix $t\in \mc S_2$. Let $n_0$  be such that $s_n\cap t = s\cap t$ for every $n\ge n_0$. Then $\langle s_n, t\rangle= \langle s,t\rangle$ for $n\ge n_0$.
 Similarly one proves the continuity with respect to the second variable.
\fprue

It is interesting to note that Talagrand's compactum given in \cite{T} can also be constructed as
$\Theta_1(\S_2)$  where  $\Theta_1: \mc S_2 \to C(\mc S)$ is the mapping $t\in \mc S_2 \mapsto
\Theta_1(t):=\Theta(\cdot,t)\in C(\mc S)$. Similarly, let $\Theta_0: \mc S\to C(\mc S_2)$ be the mapping $s\in \mc S \mapsto \Theta_0(s):=\Theta(s,\cdot)\in
C(\mc S_2)$. Set
$$ K_\al:= \Theta_0(\mc F_\al), \qquad L_\al:=\Theta_1(\mc G_\al),$$
for $\al\le \om$.

\prop
$K_\al$ and $L_\al$ have Cantor-Bendixson rank  $\al+1$ and $\om\cdot\al+1$, respectively.
\fprop

\prue
 Since  $\Theta_0$ and $\Theta_1$ are continuous,   it follows that $K_\al$  and $L_\al$ have Cantor-Bendixson index at most $\al+1$ and $\om\cdot\al+1$, respectively.  Since in addition, $\Theta_1$ is 1-1, it follows that $L_\al$ has rank exactly $\om\cdot\al +1$.

 On the other hand, let  $M:=\{2^n\}_{n\ge 0}$. Then the restriction of $\Theta_0$ to $\mc F_\al\rest M:=\conj{s\in \mc F_\al}{s\con M}$ is 1-1, which shows that the rank of $K_\al$ is $\al+1$.  Suppose that $s_0\neq s_1$ are elements of $\mc F_\al\rest M$, $s_0= \{2^{i_0}< \cdots <2^{i_{k-1}}< 2^{i_{k}}<\cdots<2^{i_l}\}$, and
       $s_1= \{2^{i_0}< \cdots <2^{i_{k-1}}< 2^{j_{k}}< \cdots< 2^{j_r}\}$, with $i_k\neq j_k$, say $i_k<j_k$. For each $m\le k$, let $t_m:=[2^{i_m}, 2^{i_m+1}-1]\in \mk S$.
       Then it follows that $2^{i_m}\in t_m$ for every $m\le k$ and $2^{j_k}\notin t_k$. If we set $t:= \bigcup_{m\le k} t_m\in \mc S_2$, then $(t_m)_{m\le k}$ is the canonical decomposition of
    $t$, hence $\langle s_0 ,t \rangle = \langle s_1,t\rangle +1$, so $\Theta(s_0,t)\neq \Theta(s_1,t)$.
\fprue

\defi
Let $E$ be a Banach space and let $\al\le \om$. By a  \emph{($\Theta$, $\al$)-embedding} into $E\times E^*$  we mean
mappings
$\ga_0: \mc F_\al\rest M\to E$ and  $\ga_1: \mc G_\al\rest M\to E^*$ for some infinite $M\con \N$ such that
\begin{enumerate}
\item[(a)] $\ga_0$ and $\ga_1$ are continuous, when we consider $E$ with its weak-topology and $E^*$ with its $\mr{weak}^*$-topology, respectively.
 \item[(b)] $\ga_1(t)\big(\ga_0(s)\big)=\Theta(s,t)$ for every $(s,t)\in \mc F_\al\rest M\times \mc G_\al\rest M$.
\end{enumerate}
The diameter of a $(\Theta,\al)$-embedding is
$$
D(\ga_0,\ga_1)=\sup\{\|\ga_0(s)\|\|\ga_1(t)\|:(s,t)\in \mc F_\al\rest M\times \mc G_\al\rest M\}.
$$
\fdefi

\teor
No reflexive Banach lattice admits a $(\Theta,\om)$-embedding. In fact, if $(\ga_0,\ga_1)$ is a $(\Theta,n)$-embedding in a reflexive Banach lattice, then $D(\ga_0,\ga_1)\ge n$.
\fteor

\prue
Suppose otherwise, and fix a $(\Theta,n)$-embedding $\ga_0: \mc F_n\rest M \to E$ and $\ga_1: \mc G_n\rest M\to E^*$. To simplify the notation, we assume without loss of generality that $M=\N$. For $s\in \mc F_n$ let $x_s=\ga_0(s)\in E$, and for $t\in \mc G_n$ let $x^*_t=\ga_1(t)\in E^*$.

\clam
There is a mapping $\De: [\N]^{\le n}\to \mk S\cup \{\buit\}$ such that
\begin{enumerate}
\item[(a)]  $\De(\buit)=\buit$, $\De(s)\neq \buit$ for every $ \#s<n$,  $s\neq \buit$, and  for every $s<m_0<m_1$,
\begin{equation*}
 n,\De(s)<\De(s\cup \{m_0\})<\De(s\cup\{m_1\}).
\end{equation*}
\item[(b)]  Let $y_\buit=x_\buit$ and $y^*_\buit=x_\buit^*$, and for each $s:=\{m_1<\cdots <m_k\}$, let
\begin{align*}
s_i:= &  \{   m_1,\dots,  m_i  \}  \text{ for $i=1,\dots, k$},  \\
y_s:= &  \frac{1}{\prod_{i=1}^k  \#(\De(s_i)  ) }  \sum_{(r_i)_{i=1}^k \in\prod_{i=1}^k    \De(s_i)    }  x_{\{r_i\}_{i=1}^k} \\
y^*_s:= & x^*_{\De(s_1)\cup\cdots \cup \De(s_k)}.
\end{align*}
\end{enumerate}
Then $y_s\in \mr{conv}\big(\ga_0(\mc F_n)\big)$ and $y^*_s\in\ga_1(\mc G_n)$, and for every $s\in \mc F_n$ with $\#s<n$ and for every $m>s$,
\begin{align}
&y_{s\cup\{k\}}\overset{\nrm{\cdot}_{L_1}}{\to}_k   y_{s}.\\
 &|y^*_{s\cup \{m\} } (y_{s}- y_{s\cup \{m\}})| = 1.
 \label{loiewrhhewwe}\\
&y^*_\buit (y_\buit)=1.
\end{align}

\fclam
\prucl
We define $\De_k: [\N]^{\le k}\to \mk S\cup \{\buit\}$ with the properties (a), (b) above by induction on $k\le n$.  If $k=0$, $\De_0(\buit)=\buit$. Suppose done for $k$, and let us define $\De_{k+1}:[\N]^{\le k+1}\to \mk S\cup \{\buit\}$ extending $\De_k$: Fix $s:=\{m_1<\cdots <m_k \}$, $s_i:=\{m_1,\dots, m_i\}$, $i=1,\dots,k$.  Then by definition,
\begin{equation}
y_s=   \frac{1}{\prod_{i=1}^k  \#(\De(s_i)  ) }  \sum_{(r_i)_{i=1}^k \in\prod_{i=1}^k    \De(s_i)    }  x_{\{r_i\}_{i=1}^k}.
\end{equation}
We know that for every  $r_1<\cdots <r_k$ one has that $x_{\{r_1,\dots,r_k, r\}}  \overset{\mr{weak},
L_1}\to_r x_{\{r_1,\dots, r_k\}}$. So by the weak-Banach-Saks property of $L_1$ \cite{Szlenk}, we can find
for every $m>m_k$ an element $\bar{s}_m\in \mk S$ with $\bar{s}_m>n$ such that
\begin{enumerate}
\item[(i)]    $\De(s)<\bar{s}_p<\bar{s}_q$ for $s<p<q$.
 \item[(ii)] For every   $(r_i)_{i=1}^k
 \in\prod_{i=1}^k    \De(s_i)$, $u:=\{r_i\}_{i=1}^k$,  one has that
 \begin{equation}
 \frac{1}{\# \bar{s}_m} \sum_{r\in \bar{s}_m} x_{u\cup \{r\}} \overset{\nrm{\cdot}_{L_1}}\to_m  x_{u}
 \end{equation}
\end{enumerate}
Let $\De(s\cup\{m\}):=\bar{s}_m$ for every $s<m$. Then it follows from (ii) that
\begin{align*}
y_{s\cup \{m\}}= & \frac{1}{\prod_{i=1}^k    \De(s_i)} \sum_{(r_i)_{i=1}^k
 \in\prod_{i=1}^k    \De(s_i)} \frac{1}{\# \bar{s}_m} \sum_{r\in \bar{s}_m} x_{\{r_i\}_{i=1}^k\cup \{r\}}  \overset{\nrm{\cdot}_{L_1}}\to_m  \\
& \overset{\nrm{\cdot}_{L_1}}\to_m   \frac{1}{\prod_{i=1}^k    \De(s_i)}  \sum_{(r_i)_{i=1}^k
 \in\prod_{i=1}^k    \De(s_i)}  x_{\{r_i\}_{i=1}^k} =y_s.
\end{align*}
For \eqref{loiewrhhewwe}:   Let $s<m$, $s=\{m_1<\ldots<m_k\}$ with $k<n$. Let $s_i=\{m_1,\ldots,m_i\}$ for $1\leq i\leq k$, and set
\begin{align*}
L_0:= &\prod_{i=1}^k \#(\De(s_i)), \\
L_1:= &  \#(\De(s\cup \{m\}),\\
L_2:= &L_0 \cdot  L_1\\
t_0:= & \bigcup_{i=1}^k \De(s_i)\\
t_1:= & t_0 \cup \De(s\cup\{m\}).
\end{align*}
Then
\begin{align*}
|y^*_{s\cup \{m\}} &(y_{s}- y_{ s\cup \{m\}})|=| x^*_{t_1} ( y_{s}- y_{s\cup \{m\}})| =   \\
=& \frac{1}{L_2}\left| \left(L_1 \cdot \sum_{v \in\prod_{i=1}^k \De(s_i)}  \ga_1(t_1)\Big(\ga_0(v)\Big) -
\sum_{w \in\big(\prod_{i=1}^k \De(s_i)\big)\times \De(s\cup \{m\})} \ga_1(t_1)\Big(\ga_0(w)\Big) \right)\right|.
\end{align*}
Now, given $v\in\prod_{i=1}^k \De(s_i)$,
\begin{equation*}
 \ga_1(t_1)\big(\ga_0(v)\big)=  \Theta(v,t_1)=\Theta(v,t_0)=\left\{\begin{array}{ll}  1 & \text{ if $k$ is even}\\
0 & \text{ if $k$ is odd.}
\end{array}\right.
\end{equation*}
Similarly, given $w \in\big(\prod_{i=1}^k \De(s_i)\big)\times \De(s\cup \{m\})$,
\begin{equation*}
 \langle \ga_1(t_1)\big(\ga_0(w)\big)=\Theta(w,t_1)= \left\{\begin{array}{ll}  0 & \text{ if $k$ is even}\\
1 & \text{ if $k$ is odd.}
\end{array}\right.
\end{equation*}
Therefore, we have
\begin{equation*}
L_1 \cdot \sum_{v \in\prod_{i=1}^k \De(s_i)}  \ga_1(t_1)\big( \ga_0(v)\big)-
\sum_{w \in\big(\prod_{i=1}^k \De(s_i)\big)\times \De(s\cup \{m\})} \ga_1(t_1)\big(\ga_0(w)\big) = \left\{\begin{array}{ll}  L_1\cdot L_0 & \text{ if $k$ is even}\\
-L_2 & \text{ if $k$ is odd.}
\end{array}\right.
\end{equation*}
Hence,
\begin{equation*}
|y^*_{s\cup \{m\}} (y_{s}- y_{ s\cup \{m\}})|=1.
\end{equation*}
\fprucl

Let $\mc U$ be a non-principal ultrafilter on $\N$. If $P(n)$ is a statement about $n\in\mathbb N$, we will write $\mc U n P(n)$ to denote that $\conj{n\in\N}{P(n)\text{ holds}}\in\mc U$. The following is a slight improvement of \cite[Proposition 1]{T}.

\clam\label{lmfkweerew3}
Suppose that $\{h_n\}$ is a bounded sequence in $E^*$. Then for every $\vep>0$ there is $a>0$ such that for every $x\in B_E$ one has that
\begin{equation}
\mc U n \int_{|h_n|>a} |h_n| |x| <\vep.
\end{equation}
\fclam
\prucl
Suppose otherwise that there is some $\vep>0$ such that for every $k$ there is $x_k\in B_E$ with
\begin{equation}
\mc U n \int_{|h_n|>k} |h_n| |x_k| >\vep.
\end{equation}
Let $L:=\sup_n \nrm{h_n}_{E^*}$. Let $y$ be a weak-accumulation point of  $\{|x_k|\}_k$.  Since $E$ is order continuous we can find $d>0$ such that $\nrm{y- y\we d}_E <\vep/(3L)$; using the equi-integrability of $B_{E^*}$, there is $a>0$ such that
\begin{equation}
\sup_n \int_{|h_n|>a } |h_n| <\frac{\vep}{3d}.
\end{equation}
Let  $a<k_1<\cdots <k_r$, and let $ (b_i)_{i=1}^r$ be a convex combination such that
\begin{equation}
\nrm{y- \sum_{i=1}^r b_i |x_{k_i}|}_E <\frac{\vep}{3L}.
\end{equation}
Choose now $n$ such that
\begin{equation}\label{oijirjwe}
\min_{1\le i\le r}\int_{|h_n|>k_i} |h_n| |x_{k_i}| >\vep.
\end{equation}
Then,
\begin{align*}
\min_{1\le i\le r}\int_{|h_n|>k_i} |h_n| |x_{k_i}|\le & \int_{|h_n|> a} |h_n| (  \sum_{i=1}^r b_i |x_{k_i}|)  \le   \int_{|h_n|> a} |h_n| y + \int |h_n|| y- \sum_{i=1}^r b_i |x_{k_i}|| \le \\
\le &   \int_{|h_n|>a} |h_n| (y\we d) + \int |h_n| |y- y\we d| + \frac{\vep}3 \le \\
\le &  d\frac{\vep}{3d}+ \nrm{h_n}_{E^*} \nrm{y-y\we d}_E+  \frac{\vep}3=\vep ,
\end{align*}
a contradiction with \eqref{oijirjwe}.
\fprucl

\clam
For every $\vep>0$ and every $0\le i \le n$ there is $s \in [\N]^i$ such that
\begin{equation}
\int |y^*_{s}| |y_{s}| \ge  1 +i -\vep(1+5i).
\end{equation}
\fclam

\prucl
Fix $\vep>0$. Induction on $0\le i\le n$. $i=0$: we know that
\begin{equation}
1-\vep<1=y^*_\buit (y_\buit)\le \int |y^*_\buit||y_\buit|.
\end{equation}
Suppose now that $s\in [\N]^i$ is such that
\begin{equation}
\int |y^*_{s}| |y_{s}| >  1 +i -\vep(1+5i).
\end{equation}
We use the previous Claim  \ref{lmfkweerew3} to find   $a>0$   such that
\begin{equation}\label{1}
\mc U m\int_{|y^*_{s\cup\{m\}}|>a} |y^*_{s\cup\{m\}}||y_s| <\vep.
\end{equation}
Since $y^*_{s\cup\{m\}}\overset{\mr{weak}}\to_my^*_s$, it follows that
\begin{equation}\label{2}
\mc U m \int|y^*_{s\cup\{m\}}||y_s|>1 +i -\vep(2+5i).
\end{equation}
Combining \eqref{1} and \eqref{2} we obtain that
\begin{equation}\label{3}
\mc U m \int_{|y^*_{s\cup\{m\}}|\le a}   |y^*_{s\cup\{m\}}||y_s|> 1 +i -\vep(3+5i).
\end{equation}
Since $y_{s\cup\{m\}} \overset{\nrm{\cdot}_{L_1}}{\to}_m y_s$, it follows from \eqref{3} that
\begin{equation}\label{4}
\mc U m \int_{|y^*_{s\cup\{m\}}|\le a}   |y^*_{s\cup\{m\}}||y_{s\cup\{m\}}|> 1 +i -\vep(4+5i).
\end{equation}
On the other hand, given $n$ we have  from \eqref{loiewrhhewwe} that
\begin{align}
\int_{|y^*_{s\cup\{m\}}|>a} |y^*_{s\cup \{m\}}||y_{s\cup \{m\}}| \ge &  \int |y^*_{s\cup\{m\}}| |y_s-y_{s\cup\{m\}} | - \int_{|y^*_{s\cup\{m\}}|\le a}|y^*_{s\cup\{m\}}| |y_s  -y_{s\cup\{m\}}| - \nonumber \\
- &  \int_{|y^*_{s\cup\{m\}}|>a} |y^*_{s\cup\{m\}}| |y_s|\ge \nonumber \\
\ge &1 -   \int_{|y^*_{s\cup\{m\}}|\le a}|y^*_{s\cup \{m\}}| |y_s  -y_{s\cup \{m\}}| - \int_{|y^*_{s\cup \{m\}}|>a} |y^*_{s\cup \{m\}}| |y_s |. \label{11}
\end{align}
Since $y_{s\cup \{m\}}\overset{\nrm{\cdot}_{L_1}}{\to}_m y_s$, it follows that
\begin{equation}\label{12}
 \mc U m \int_{|y^*_{s\cup \{m\}}|\le a}|y^*_{s\cup \{m\}}| |y_s  -y_{s\cup \{m\}}| <\vep.
\end{equation}
Combining \eqref{1} and \eqref{12} in \eqref{11} we obtain that
\begin{equation}\label{13}
\mc U m \int_{|y^*_{s\cup\{m\}}|>a} |y^*_{s\cup \{m\}}| |y_{s\cup\{m\}}|>1-2\vep.
\end{equation}
Combining \eqref{4} and \eqref{13}, we obtain
\begin{equation}\label{14}
 \mc Um \int |y^*_{s\cup \{m\}}| |y_{s\cup \{m\}}|>1 +i +1 -\vep(1+5(i+1)).
\end{equation}

\fprucl

From this last claim, we can find $s\in [\N]^n$ such that
\begin{equation}
\int |y^*_s y_s|\ge n.
\end{equation}
Since, $y_s\in\mr{conv}\big(\ga_0(\mc F_n)\big)$ and $y_s^*\in \ga_1(\mc G_n)$, it follows that $D(\ga_0,\ga_1)\geq n$.
\fprue

\prop
Let $E$ be a reflexive Banach space. \begin{enumerate}
\item[(a)] Suppose $T:E\rightarrow C(\mc S_2)$ is such that $K_\al\con T(B_E)$. Then there is a $(\Theta,\al)$-embedding in $E\times E^*$ with diameter less or equal than $\|T\|$.
\item[(b)] Suppose $T:E\rightarrow C(\mc S)$ is such that $L_\al\con T(B_E)$. Then there is a $(\Theta,\al)$-embedding in $E^*\times E$ with diameter less or equal than $\|T\|$.
 \end{enumerate}
 \fprop
\prue
Suppose that $T: E\to C(\mc S_2)$ is such that $K_\al\con T(B_E)$. For each maximal $s\in \mc F_\al$ choose $x_s\in B_E$ such that $T(x_s)= \Theta_0(s)$. From Lemma \ref{weak tree} there is an infinite $M\con \N$  such that $(\ga_0(s))_{s\in \mc F_\al\rest M}$ is a weakly convergent tree.  Since $T$ is weak-weak-continuous, it follows that $T(\ga_0(s))=\Theta_0(s)$ for every $s\in \mc F_\al \rest M$.     For each $t\in \mc G_\al$, let  $\ga_1(t):= T^*(\de_t)$.    Then for every $(s,t)\in \mc F_\al\rest M\times \mc G_\al\rest M$ one has that
\begin{equation}
\ga_1(t)\big(\ga_0(s)\big)  =  T^*(\de_t) \big(\ga_0(s)\big)= \de_t \big(T(\ga_0(s))\big)= \de_t\big( \Theta_0(s)\big) = \Theta(s,t).
\end{equation}
Notice that $\|\ga_0(s)\|\|\ga_1(t)\|\leq\|T^*\|=\|T\|$, thus $D(\ga_0,\ga_1)\leq \|T\|$.

Now, if $T: E\to C(\mc S)$ is such that $L_\al\con T(B_E)$, then arguing as above, one could find $M\con \N$ and  a continuous $\ga_1: \mc G_\al\rest M\to  E$ such that $T(\ga_1(t))=\Theta_1(t)$ for every $t\in \mc G_\al\rest M$. Let $\ga_0: \mc G_\al \rest M\to E^*$, $\ga_0(s):=T^*(\de_s)$ for every $s\in \mc G_\al\rest M$.  Again, $ \ga_1(t)\big( \ga_0(s)\big)  =  \Theta(s,t)$, and we get $D(\ga_0,\ga_1)\leq \|T\|$.
\fprue

As a consequence we get the following:

\cor
Neither $K_\om$ nor $L_\om$ are shellable by reflexive Banach lattices. In fact, if $T:E\rightarrow C[0,1]$ is such that $K_n$ or $L_n\subseteq T(B_E)$, then $\|T\|\ge n$. \qed
\fcor

Several important aspects of the geometric structure of $C(\Omega)$ spaces are closely related to operators $T:C(\Omega)\to C(\Omega)$ (see \cite{Gas, Ros} for recent and complete accounts on this relation). Note that the compact $K_\omega$ and $L_\om$ constructed above are not contained in the image of the unit ball by any weakly compact operator $T:C[0,1]\rightarrow C[0,1]$. In fact, we have the following:

\prop
Let $\Omega$ be a compact Hausdorff space and $K\subset C(\Omega)$ a weakly compact set such that $K\subseteq T(B_{C(\Omega)})$ for some
weakly compact operator $T:C(\Omega)\rightarrow C(\Omega)$. Then $K$ is in
$\mr{Sh}_{\mc{RL}}(C(\Omega))$.
\fprop

\prue
Since $T$ is weakly compact, so is its adjoint. Now, $C(\Omega)^*$ does not contain $c_0$ so the solid hull of
$T(B_{C(\Omega)^*})$ is weakly compact. Thus, we could factor $T^*$, and also $T$, through a reflexive Banach
lattice.
\fprue

The converse however is not true. In fact, it is not even true that a subset of $C[0,1]$ shellable by a
reflexive lattice is contained in $T(B_{C[0,1]})$ for a \emph{weak Banach-Saks} operator $T:C[0,1]\to
C[0,1]$. Recall that an operator is  weak Banach-Saks  when it  sends weakly compact sets to Banach-Saks
sets. Note that the Dunford-Pettis property of $C[0,1]$ yields in particular that every weakly compact operator $T:C[0,1]\to
C[0,1]$ is weak Banach-Saks. Now, suppose that   $X$ is a copy of the Schreier space in $C[0,1]$ and let $K=\{0\}\cup\{u_n\}_n$ where
$(u_n)_n$ is the unit basis of $X$. By Theorem \ref{kfksdfsklfsd}, we have that $K$ is in $\mr{Sh}_{\mc{RL}}(C[0,1])$.
If $T:C[0,1]\to C[0,1]$ is such that $K\con T(B_{C[0,1]})$, then one can prove that its  \emph{Szlenk index}
(see \cite{Ros}) is at least $\om^2$. Hence, by a result of D. E. Alspach \cite{Als}, $T$ fixes a copy of
$C(\mc S)$, hence $T$ cannot be weak Banach-Saks.

It was observed in \cite{FT} that Talagrand's weakly compact set $K_{\mc T}$ is a Banach-Saks set. Since
every space with the Banach-Saks property is reflexive, as a consequence of \cite{T} it follows that $K_{\mc
T}$ is not shellable by Banach lattices with the Banach-Saks property. In a similar spirit, one might wonder
what is the smallest ordinal $\alpha$ such that there exists a Banach-Saks set $K\subseteq C[0,1]$
homeomorphic to $\alpha$ which is not shellable by Banach lattices with the Banach-Saks property. It can be
seen that the compact set $K_\omega$ constructed before fails the Banach-Saks property. However, in
\cite{LRT} an example was given of a Banach-Saks set whose convex hull is not Banach-Saks, so it is not even
shellable by Banach spaces with the Banach-Saks property. Note this set is homeomorphic to $\omega+1$.

\end{document}